\documentclass[aps,pra,a4paper,10pt,notitlepage]{revtex4-2}

\usepackage[utf8]{inputenc}
\usepackage[T1]{fontenc}

\usepackage{appendix}
\usepackage{graphicx}
\usepackage{amsmath}
\usepackage{amsfonts}
\usepackage{amssymb}
\usepackage{xcolor}
\usepackage{fdsymbol}
\usepackage{hyperref}
\usepackage{pifont}
\newcommand{\Ph}{\phantom{-}}
\graphicspath{{./figures/}}

\newcommand{\revision}{\color{black}}

\definecolor{O}{HTML}{d32424}
\definecolor{C}{HTML}{24d624}
\definecolor{H}{HTML}{cacaca}

\begin{document}
\title{Roaming in acetaldehyde}

 \begin{abstract}
We investigate roaming in the photodissociation of acetaldehyde (CH$_3$CHO), providing insight into the contrasting roaming dynamics observed for this molecule compared to formaldehyde.  We carry out trajectory studies for full-dimensional acetaldehyde, supplemented with an analysis of a two degree-of-freedom restricted model and obtain 
evidence for two distinct roaming pathways. Trajectories exhibit roaming at both shorter (9-11.5 au) and larger (14.5-22.9 au) maximum CH$_3$-HCO separations, characterized by differing amounts of HCO rotation. No roaming trajectories were found in the intervening gap region. The roaming dynamics near 14.5-22.9 au are well-reproduced by the restricted model and involve passage through a centrifugal barrier, analogous to formaldehyde roaming.  However, the shorter-range 9-11.5 au roaming appears unique to acetaldehyde, and is likely facilitated by repulsive interactions absent in the simplified models. Phase space analysis reveals that  this additional roaming pathway is inaccessible in the reduced dimensionality system. The findings suggest acetaldehyde's increased propensity for roaming compared to formaldehyde may arise from the presence of multiple distinct roaming mechanisms rather than solely the higher roaming fragment mass.

 \end{abstract}

 \author{Vladim{\'i}r Kraj{\v{n}}{\'a}k}
 \email{Corresponding author; vkrajnak@gmail.com}
 \address{School of Mathematics, University of Bristol, Fry Building, Woodland Road, Bristol, BS8 1UG, United Kingdom.}

 \author{Stephen Wiggins}
 \address{School of Mathematics, University of Bristol, Fry Building, Woodland Road, Bristol, BS8 1UG, United Kingdom.\\ 
Department of Mathematics, 
United States Naval Academy, Chauvenet Hall, 572C Holloway Road,
Annapolis, MD 21402-5002, USA.}

\maketitle

\thispagestyle{empty}

\section{Introduction}
\label{sec:intro}

 Anomalies in experimental data for the photodissociation of formaldehyde, ${\rm H_2CO}$, \cite{vanzee1993} have led to new understanding of the dissociation dynamics of formaldehyde \cite{townsend2004roaming, lahankar2006roaming, lahankar2007energy, lahankar2008further, shepler2011roaming, houston2016roaming} which has been referred to as the {\em roaming mechanism} for reaction dynamics. The widespread interest in this new reaction mechanism resulted  in its discovery in a variety of reactions. Much of this work was described  in the abundance of review papers that subsequently appeared \cite{suits2008roaming, bowmansuits2011roaming, bowmanshepler2011roaming, bowman2014roaming, bowman2017theories, mauguiere2017roaming, suits2020roaming}.

 This begs the question, ``what is roaming''?  In \cite{mauguiere2014multiple} the authors noted that a dissociating molecule should possess two essential characteristics in order to label the reaction as ``roaming''. In particular, the molecule should have competing dissociation channels, such as dissociation to molecular and radical products, and there should exist a long range attraction between fragments of the molecule. In a recent review article \cite{suits2020roaming} Suits refines this definition of roaming even further by stating that ``A roaming reaction is one that yields products via reorientational motion in the long-range ($3-8$\AA) region of the potential''. 

 Following the seminal studies on roaming in the dissociation of formaldehyde, attention was focused on roaming in acetaldehyde, ${\rm CH_3 CHO }$\cite{houston2006photodissociation, shepler2007, shepler2008roaming, heazlewood2008roaming, lee2014two, li2015communication, rubio2012imaging, rubio2007slice, Harding2010roaming, Klippenstein2011roaming}. Formaldehyde and acetaldehyde roaming dynamics present interesting  contrasts. In formaldehyde the long range reorientational motion occurs between the HCO fragment and a hydrogen atom. In acetaldehyde the long range reorientation dynamics occurs between the HCO fragment and a methyl group, ${\rm CH_3}$. Both molecules possess competing dissociation channels to molecular and radical products. It was found that  formaldehyde dissociation by roaming occurred much less frequently than dissociation to molecular products over a potential saddle point. In contrast, acetaldehyde prefers to dissociate to molecular products via roaming.
 
 A reason for this could be the mass difference in the roaming fragment; H for formaldehyde and ${\rm CH_3}$ for acetaldehyde. This question was investigated in \cite{krajnak2018influence} using a model that has been shown to exhibit the essential features of roaming, the Chesnavich model \cite{chesnavich1986multiple} for the reaction CH$_4^+ \rightarrow$ CH$_3^+$ + H  \cite{mauguiere2014multiple, mauguiere2014roaming, krajnak2018phase}. An advantage of the Chesnavich model is that it is analytical and  therefore the parameters that describe the reaction can be varied. In particular, the mass of the roaming fragment, H can be varied and we can study the effect of this variation  on the roaming mechanism to dissociation.  Following a detailed analysis of trajectories and phase space structures responsible for roaming, it was found that the variation of the mass of the roaming fragment had no significant effect on the quantity of trajectories dissociating via the roaming channel. Perhaps the simplicity of Chesnavich's model is preventing a comparison of features mimicking formaldehyde and acetaldehyde. Recent work \cite{yang2020real} raises questions about the nature and quantity of trajectories undergoing dissociation of acetaldehyde due to roaming at certain energies.

 A geometrical phase space approach for revealing the roaming mechanism was previously carried out for a reduced dimensionality formaldehyde model \cite{mauguiere2015phase} and its validity was verified in full-dimensional studies of formaldehyde \cite{houston2016roaming}. In that work it was shown that the stable and unstable manifolds of certain periodic orbits determined whether or not a trajectory dissociated to molecular or radical products. Due to the unboundedness of the energy surface that is characteristic for all dissociation problems, the ergodic hypothesis fails for roaming systems and invalidates standard statistical methods for reaction dynamics.

 In this work we investigate roaming in acetaldehyde. Due to high dimensionality, we perform the geometrical phase space analysis on a highly restricted model of acetaldehyde with 2 degrees of freedom and validate the findings in the full-dimensional model of acetaldehyde. We find roaming trajectories that are in agreement with the findings of \cite{shepler2007, shepler2008roaming, Harding2010roaming}, yet their properties suggest two distinct roaming pathways.

 This paper is outlined as follows. In Section \ref{sec:charac} we describe the nature of the roaming trajectories that we have computed. In Section \ref{sec:methods} we describe our method for determining and characterizing roaming trajectories. In Section \ref{sec:pert}  we describe our computational approach for identifying roaming trajectories and in Section \ref{sec:compare} we compare our full dimensional calculations with calculations for a two degree-of-freedom restricted system. In Section \ref{sec:concl} we describe our conclusions for the full dimensional trajectory analysis of roaming. In Appendix \ref{sec:2dofmodel} we describe our two degree-of-freedom model of restricted acetaldehyde. In Section \ref{sec:identify} we show how to identify the relevant periodic orbits governing roaming and discuss their geometry. In Section \ref{sec:rra} we discuss the  ''phase space template'' governing the roaming process, which includes finding the relevant periodic orbits (Section \ref{sec:finding}) and then using them to construct the relevant phase space dividing surfaces (Section \ref{sec:DS}). In Section \ref{sec:concl_res} we give our conclusions for roaming in restricted acetaldehyde.

\section{Characterization of Roaming Trajectories}
\label{sec:charac}

In this section we describe the roaming  trajectories that we have found
from our  investigation of a full-dimensional acetaldehyde model at zero total angular momentum and an energy $0.1612$ a.u. above the bottom of the acetaldehyde well at $-153.6017$ a.u. We use the potential energy surface used in \cite{Han2017} and based on \cite{heazlewood2008roaming} kindly provided by Prof. J. Bowman, Prof. Y.-C. Han and Prof. B. Fu. Of interest are trajectories where CH$_3$ roams in a region of nearly constant potential energy and abstracts H from HCO prior to dissociation. For each (discretised) roaming trajectory we record the largest distance between the centres of mass of CH$_3$ and HCO. We find two disjoint intervals with approximate endpoints $9-11.5$ a.u. and $14.5-22.9$ a.u., where for each value in the interior of the two intervals we can find a roaming trajectory which reaches that maximum roaming distance. We find no roaming trajectories between $11.5$ and $14.5$ a.u.
 
 All roaming trajectories are characterised by a larger C-C distance at the moment when H is transferred \cite{shepler2007}, than trajectories passing near the potential saddle point, although we do not have enough data for statistically significant conclusions and quantitative results are outside of the scope of this work.

\begin{figure}
\centering
 \vspace*{-5mm}
 \includegraphics[width=\textwidth]{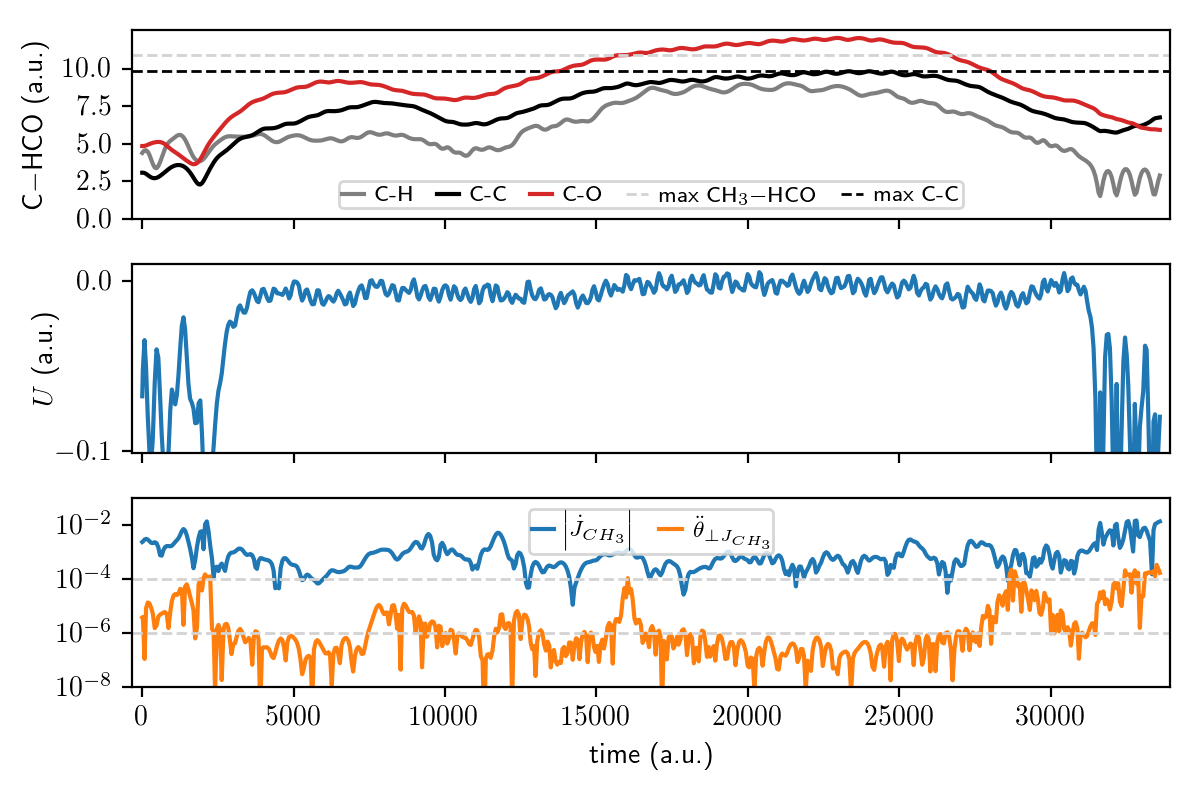}
 \vspace*{-5mm}
 \caption{Roaming trajectory, short CH$_3$-HCO distance, O (red) always faces away from CH$_3$. {\revision The upper panel displays the evolution of distances between C of the CH$_3$ fragment and H, C, O; as well as the maximal C-C distance and maximal distance between the centres of mass of CH$_3$ and HCO. The middle panel displays the potential energy along this trajectory. The bottom panel displays ($y$ axis in logarithmic scale) the time derivatives of the unit vector $J_{CH_3}$ and the angle $\theta_{\perp J_{CH_3}}$ defined in Sec. \ref{sec:compare}, which indicate out-of-plane motion and the coupling between the angular and radial degrees of freedom.} Trajectory data available in Supplementary material, animation published under https://youtu.be/mtpSFWB7Tug.}
 \label{fig:roaming_short_O_away}
\end{figure}

\begin{figure}
\centering
 \vspace*{-5mm}
 \includegraphics[width=\textwidth]{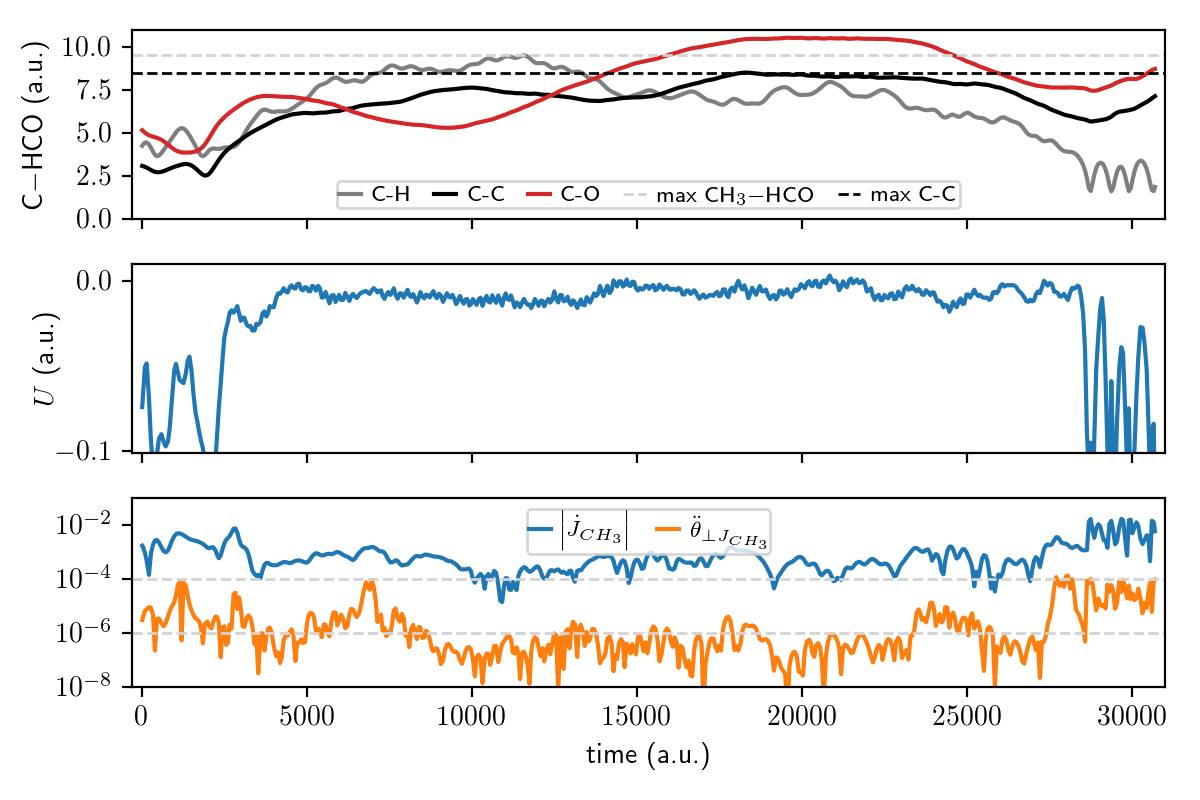}
 \vspace*{-5mm}
 \caption{Roaming trajectory, short CH$_3$-HCO distance, O (red) faces CH$_3$ between $t=0$ and $t=5000$. {\revision The upper panel displays the evolution of distances between C of the CH$_3$ fragment and H, C, O; as well as the maximal C-C distance and maximal distance between the centres of mass of CH$_3$ and HCO. The middle panel displays the potential energy along this trajectory. The bottom panel displays ($y$ axis in logarithmic scale) the time derivatives of the unit vector $J_{CH_3}$ and the angle $\theta_{\perp J_{CH_3}}$ defined in Sec. \ref{sec:compare}, which indicate out-of-plane motion and the coupling between the angular and radial degrees of freedom.} Trajectory data available in Supplementary material, animation published under https://youtu.be/fgTzZ-avjtY.}
 \label{fig:roaming_short}
\end{figure}

\begin{figure}
\centering
 \includegraphics[width=\textwidth]{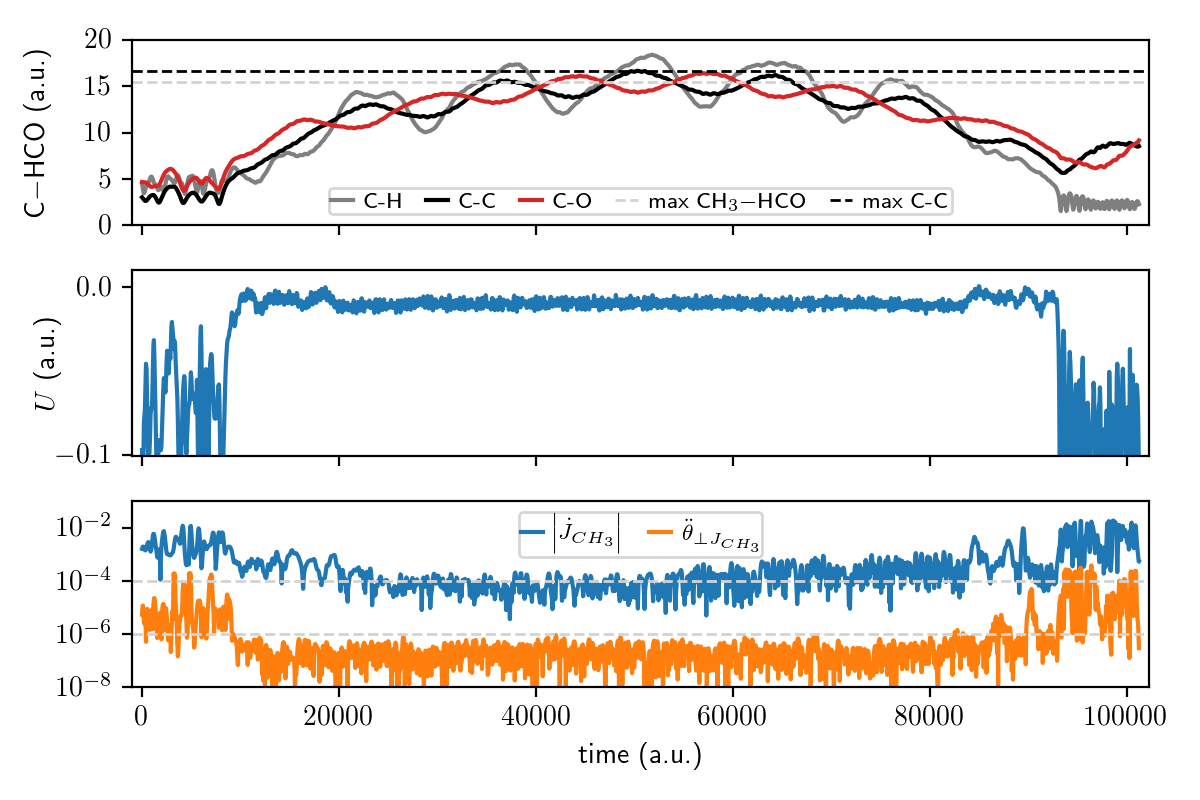}
 \vspace*{-5mm}
 \caption{Roaming trajectory, large CH$_3$-HCO distance. {\revision The upper panel displays the evolution of distances between C of the CH$_3$ fragment and H, C, O; as well as the maximal C-C distance and maximal distance between the centres of mass of CH$_3$ and HCO. The middle panel displays the potential energy along this trajectory. The bottom panel displays ($y$ axis in logarithmic scale) the time derivatives of the unit vector $J_{CH_3}$ and the angle $\theta_{\perp J_{CH_3}}$ defined in Sec. \ref{sec:compare}, which indicate out-of-plane motion and the coupling between the angular and radial degrees of freedom.} Trajectory data available in Supplementary material, animation published under https://youtu.be/AxMpUb6xLxE.}
 \label{fig:roaming_15}
\end{figure}

\begin{figure}
\centering
 \includegraphics[width=\textwidth]{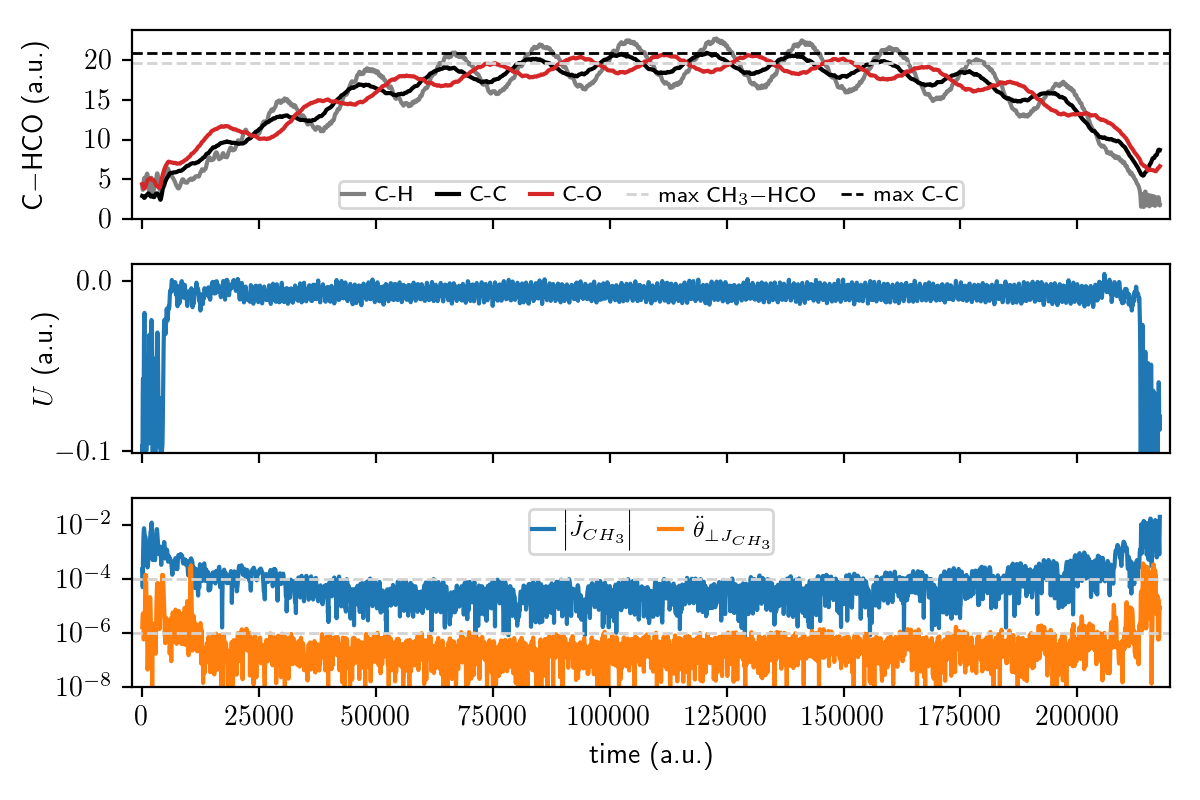}
 \vspace*{-5mm}
 \caption{Roaming trajectory passing near the centrifugal barrier, which is located at around 22.9 a.u. {\revision The upper panel displays the evolution of distances between C of the CH$_3$ fragment and H, C, O; as well as the maximal C-C distance and maximal distance between the centres of mass of CH$_3$ and HCO. The middle panel displays the potential energy along this trajectory. The bottom panel displays ($y$ axis in logarithmic scale) the time derivatives of the unit vector $J_{CH_3}$ and the angle $\theta_{\perp J_{CH_3}}$ defined in Sec. \ref{sec:compare}, which indicate out-of-plane motion and the coupling between the angular and radial degrees of freedom.} Trajectory data available in Supplementary material, animation published under https://youtu.be/VmXUT1ahZ3M.}
 \label{fig:roaming_21}
\end{figure}

 We present four roaming trajectories: two roam at a maximal distance in the $9-11.5$ a.u. interval, Fig. \ref{fig:roaming_short} and Fig. \ref{fig:roaming_short_O_away}, and two roam at a maximal distance in the $14.5-22.9$ a.u. interval, Fig. \ref{fig:roaming_15} and Fig. \ref{fig:roaming_21}. The figures show the evolution of distances between C and all HCO atoms in time, complemented by the maximal roaming distance between CH$_3$ and HCO centres of mass and potential energy along the trajectory. On the lower end of {\revision the time scale}, all trajectories have a visible C-C bond manifested by a short C-C distance and low potential energy. On the higher end, all trajectories show a C-H bond and a significant distance between C and CO. In between, CH$_3$ roams in an area of nearly constant potential (up to internal vibrations).

\begin{figure}
\centering
 \includegraphics[width=\textwidth]{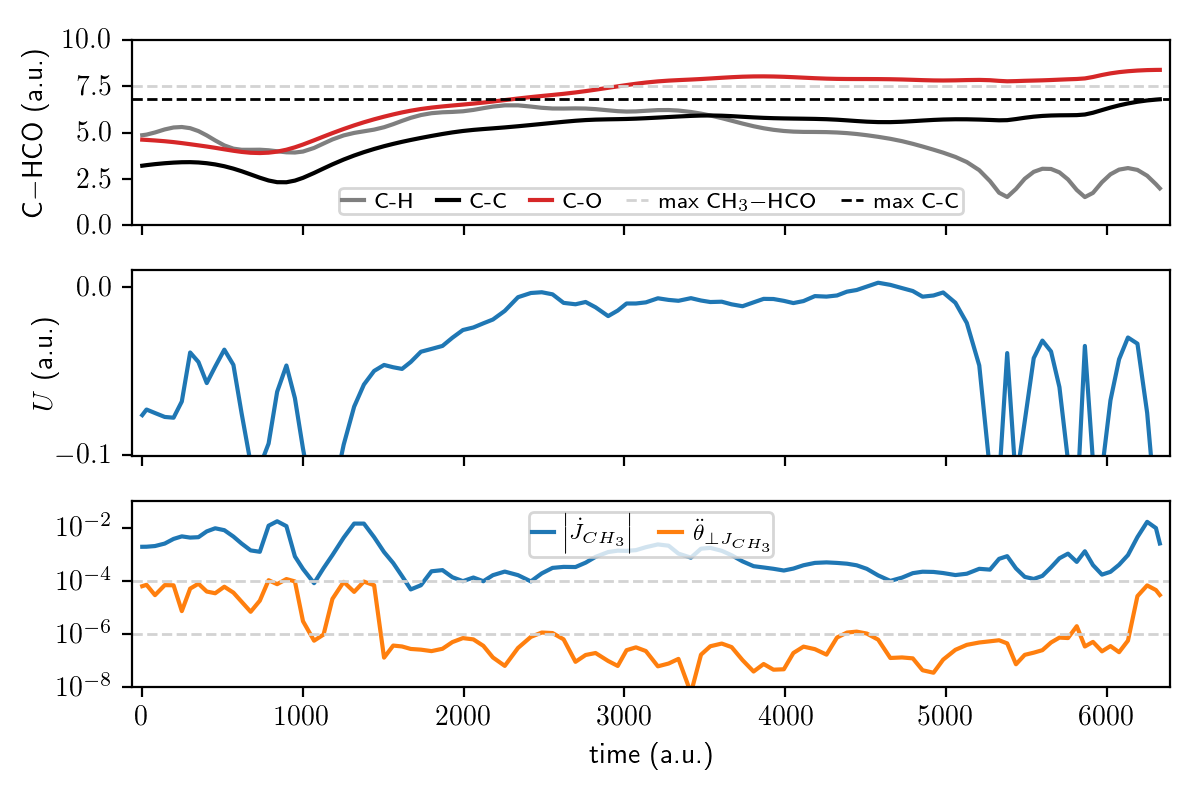}
 \caption{Dissociation trajectory passing near the saddle point. {\revision The upper panel displays the evolution of distances between C of the CH$_3$ fragment and H, C, O; as well as the maximal C-C distance and maximal distance between the centres of mass of CH$_3$ and HCO. The middle panel displays the potential energy along this trajectory. The bottom panel displays ($y$ axis in logarithmic scale) the time derivatives of the unit vector $J_{CH_3}$ and the angle $\theta_{\perp J_{CH_3}}$ defined in Sec. \ref{sec:compare}, which indicate out-of-plane motion and the coupling between the angular and radial degrees of freedom.} }
 \label{fig:saddle}
\end{figure}

 For completeness of dissociation mechanisms of acetaldehyde into CH$_4$ and CO, we include a trajectory passing near the potential saddle point in Fig. \ref{fig:saddle}.
 
 There are multiple features of the presented roaming trajectories that suggest they are due to two separate mechanisms:
 \begin{itemize}
  \item they explore different areas of the potential energy surface,
  \item different lengths of trajectories,
  \item only one kind of behaviour is captured by a restricted 2 DoF system,
  \item a gap in which we were not able to find roaming with our approach.
 \end{itemize}

 As represented by the trajectory shown in Fig. \ref{fig:roaming_short_O_away}, some roaming trajectories in the $9-11.5$ a.u. interval are characterised by a HCO in which O always faces away from CH$_3$. Several roaming trajectories we found roaming between $9-11.5$ a.u. have this property, with an exception shown in Fig. \ref{fig:roaming_short}. In contrast, HCO rotates significantly more in the C-C-O plane in trajectories that roam between $14.5-22.9$ a.u., that is near the centrifugal barrier located at approximately $22.9$ a.u. Examples are shown in Fig. \ref{fig:roaming_15} and Fig. \ref{fig:roaming_21}. 
 
 Trajectories in the $9-11.5$ a.u. and $14.5-22.9$ a.u. intervals also differ in their lengths. We found roaming trajectories as short as $14000$ time units (a.u.) between $9-11.5$ a.u., while the shortest near the centrifugal barrier is $80000$ time units long. We remark that according to dynamical systems theory, there is no upper bound on the length of trajectories that visit the neighbourhood of invariant objects such as the one due to the centrifugal barrier.
 
 We compared the dynamics in the full-dimensional acetaldehyde model to a restricted 2 DoF system, the details of which we present in Methods and Appendix. We found that the restricted system, in which we can determine all reactive mechanisms using dynamical systems theory, only exhibited roaming near the centrifugal barrier. This is in agreement with earlier studies of formaldehyde and Chesnavich's CH$_4^+$ model \cite{mauguiere2014multiple,mauguiere2014roaming, mauguiere2015phase, mauguiere2017roaming, krajnak2018phase} and was shown to be robust under variations of certain parameters \cite{krajnak2018influence}.
 
 Lastly, we were not able to find roaming in the gap between the two intervals. Instead in the gap we find trajectories that are trapped in an area of nearly constant potential for some time, but never form abstract H from HCO or never form acetaldehyde.

\section{Methods for Identifying Roaming Trajectories}
\label{sec:methods}

 The results described above rely on two separate and independent approaches that complement each other in pursuit of the understanding of roaming in acetaldehyde. On one hand, we analyse the phase space structures in a highly restricted model of acetaldehyde with 2 degrees of freedom. This allows us to relate the restricted acetaldehyde model to other models with 2 degrees of freedom, where roaming and related aspects of dynamics were explained using invariant structures in phase space \cite{mauguiere2014multiple,mauguiere2014roaming, mauguiere2015phase, mauguiere2017roaming, krajnak2018phase, krajnak2018influence, krajnak2019isokinetic, krajnak2020manifld}. On the other hand, we explore the full-dimensional acetaldehyde system by repetitively perturbing non-dissociative trajectories that initially resemble trajectories leading to radical products. This way we are able to relate observations to existing results on acetaldehyde \cite{houston2006photodissociation, shepler2007, shepler2008roaming, heazlewood2008roaming, lee2014two, li2015communication, rubio2012imaging, rubio2007slice, Harding2010roaming, Klippenstein2011roaming} and validate the low dimensional phase space analysis.
 
 Dynamics in the full-dimensional and restricted systems are closely related for large CH$_3$ - HCO distances, most notably the proof for the existence of a centrifugal barrier due to the decoupling of angular and radial degrees of freedom \cite{krajnak2018phase} applies to both. Following dynamical systems theory, a normally hyperbolic invariant manifold (torus) is associated with a centrifugal barrier and this object, via its stable and unstable invariant manifolds (cylinders \cite{Wiggins01,Uzer02}), allows certain trajectories to dissociate into radical products and prevents other trajectories from doing so. The latter is a fundamental ingredient of roaming.
  
 \subsection{Perturbative search}
 \label{sec:pert}
 
 The perturbative search of non-dissociative trajectories in the full-dimensional system relies on what is described as follows. We identified several trajectories that do not reach a CH$_3$ - HCO distance corresponding to the centrifugal barrier (approximately $22.9$ a.u.) within $10000$ time units. From these trajectories we generated new initial conditions by adding random perturbations sampled from a multivariate Gaussian distribution and integrating these forward and backward in time. By imposing stronger constraints, such as longer non-dissociation time and at least one ends of the trajectory to be an acetaldehyde configuration, we found a variety of trajectories that include roaming and dissociation over the potential saddle. All presented trajectories are related by a (long) sequence of perturbations.
 
 This approach is feasible because transport in Hamiltonian dynamical systems is driven by stable and unstable invariant manifolds (cylinders) of hyperbolic objects that delimit regions in phase space. Unstable invariant manifolds carry trajectories out of a region, while stable invariant manifolds carry trajectories to a region. Thus to break the C-C bond in acetaldehyde and form a C-H bond between CH$_3$ and HCO, all corresponding trajectories `start' and `end' enclosed by the same invariant manifolds (cylinders). Following dynamical systems theory, dynamics near the intersection of stable and unstable invariant manifolds is chaotic and thus any to trajectories that are at some point in time arbitrarily close may evolve in completely different ways. Therefore it is not surprising that by slightly perturbing a trajectory, we may change the part of phase space the trajectory explores and the mechanism that determines its evolution (radical dissociation, roaming, dissociation over a saddle point) and properties (time, bonds, internal energy redistribution).
 
 We remark that aside from trajectories exhibiting roaming and dissociating over the potential saddle, we observed trajectories that, despite never reaching the acetaldehyde well or never abstracting H from HCO, spent a considerable amount of time in the area of nearly constant potential, similar to behaviour that was refered to as transient roaming in \cite{Endo2020}.
 
 \subsection{Comparison of full-dimensional and restricted dynamics}
 \label{sec:compare}
 
 The highly restricted system with 2 degrees of freedom serves the purpose of allowing an exhaustive investigation of phase space structures in a system that is, as much as 2 degrees of freedom permit, related to acetaldehyde. It treats CH$_3$ and HCO as rigid bodies and the centre of mass of CH$_3$ moves in the HCO plane. The restricted system possesses three potential wells, one near each atom of HCO and two of which, the ones near C and H, are accessible at the total energy we used, namely $0.01$ a.u. above the restricted radical dissociation threshold at $-153.4405$ a.u. This is the same total energy as considered in the full-dimensional system.
 
 We find unstable periodic orbits, each of which is at the narrowest part of a different bottleneck on the surface of constant energy, that naturally divide the space into regions. One is at the edge of each potential well corresponding to the formation/breaking point of C-C and C-H bonds. The centrifugal barrier is located at approximately $22.9$ a.u. and marks the location beyond which acetaldehyde has dissociated to radical products. These three orbits were also identified in Chesnavich's phenomenological CH$_4^+$ model \cite{mauguiere2014multiple,mauguiere2014roaming} and a similarly restricted formaldehyde \cite{mauguiere2015phase}. In contrast to those systems, the restricted acetaldehyde has an additional unstable periodic orbit (and bottleneck) separating the two potential wells from what is usually referred to as the flat region. This bottleneck enables the reduced system to have trajectories that correspond to dissociation to molecular product over the potential saddle. Due to the absence of bottlenecks in the flat region, all roaming trajectories pass `near' the centrifugal barrier while they roam, making the roaming mechanism in the restricted system equivalent to Chesnavich's model and restricted formaldehyde. A detailed phase space analysis can be found in the Appendix.
  
 Qualitatively, the similarity of full-dimensional and restricted systems for large CH$_3$ - HCO distances can be seen from two facts in the presented figures. Firstly, the unit vector of the angular momentum of CH$_3$ relative to the centre of mass of whole molecule ($J_{CH_3}$) varies little for large CH$_3$ - HCO distances. The norm of the time derivative of $J_{CH_3}$ is $\approx 10^{-4}$, see Fig. Fig. \ref{fig:roaming_15} and Fig. \ref{fig:roaming_21}. In contrast, for shorter CH$_3$ - HCO distances this quantity is around $10^{-3}$ or $10^{-2}$. Thus near the centrifugal barrier CH$_3$ moves with respect to HCO on a nearly fixed plane.
 
 Secondly, consider the angle $\theta_{\perp J_{CH_3}}$ defined by the vectors joining C and O, and the centres of mass of HCO and CH$_3$ projected onto the plane of motion. If the angular and radial degrees of freedom were fully decoupled, this angle would evolve linearly in time, thus $\ddot\theta_{\perp J_{CH_3}} = 0$. {\revision The degrees of freedom are never fully decoupled along the presented trajectories, but $\ddot\theta_{\perp J_{CH_3}}$ stays smaller than $ 10^{-6}$ for large CH$_3$ - HCO distances in Fig. \ref{fig:roaming_15} and Fig. \ref{fig:roaming_21}. For CH$_3$ - HCO distances between $9-11.5$ a.u., $\ddot\theta_{\perp J_{CH_3}}$ is slightly larger as shown in Fig. \ref{fig:roaming_short_O_away} and Fig. \ref{fig:roaming_short}. More importantly,} $J_{CH_3}$ varies significantly more, suggesting that out-of-plane motion is important for these roaming trajectories. For short CH$_3$ - HCO distances $\ddot\theta_{\perp J_{CH_3}}$ is near $10^{-5}$ or larger.
 
 If CH$_3$ moves with respect to HCO on a nearly fixed plane and the angle $\theta_{\perp J_{CH_3}}$ on this plane evolves nearly linearly, the system's evolution near the centrifugal barrier only depends on the CH$_3$ - HCO distance. In this area of phase space the system is effectively a 1 degree-of-freedom system, plus internal vibrational degrees of freedom of CH$_3$ and HCO, which do not influence the motion of the centres of mass. Therefore it can be reasonably expected that a system with 2 degrees of freedom can capture the essential aspects of roaming in this area of phase space.
 
\begin{figure}
\centering
 \includegraphics[width=\textwidth]{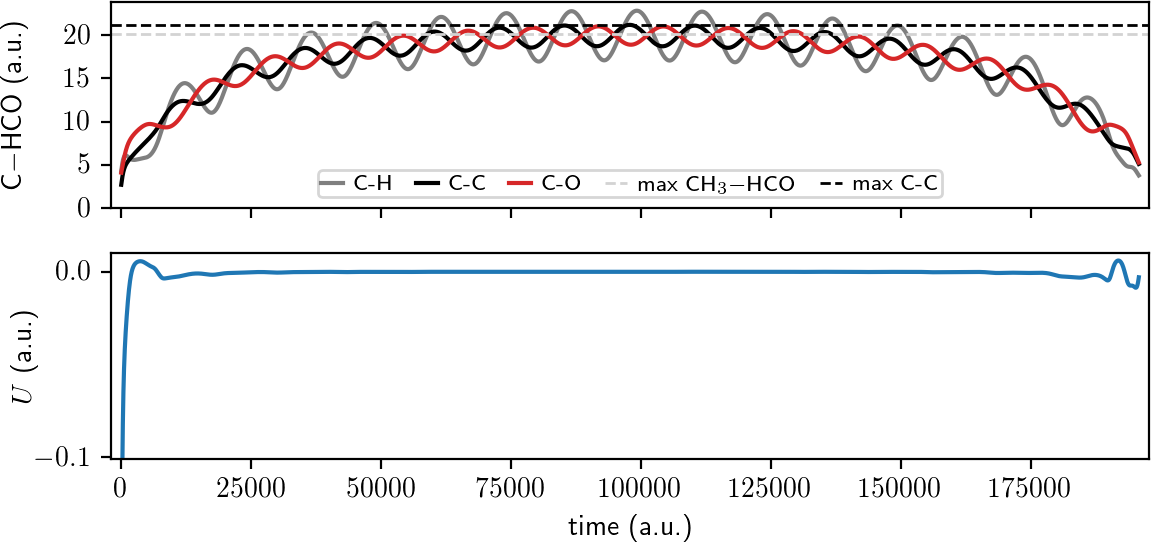}
 \vspace*{-5mm}
 \caption{Roaming trajectory of the restricted system passing near the centrifugal barrier, which is located at around 22.9 a.u. The trajectory was initialised using a point on the trajectory shown in Fig. \ref{fig:roaming_21}.}
 \label{fig:roaming_restricted_21}
\end{figure}

 To illustrate the similarity, we present a roaming trajectory of the restricted system in Fig. \ref{fig:roaming_restricted_21} that was initialised using a point on the roaming trajectory of the full-dimensional system shown in Fig. \ref{fig:roaming_21}.
  
 While roaming mechanisms in the restricted models of acetaldehyde and formaldehyde are nearly identical, experimental and numerical investigations of their full-dimensional counterparts tell very different stories. This fact and the above-presented trajectories of full-dimensional acetaldehyde suggest that formaldehyde's roaming mechanism is present in acetaldehyde, but another mechanism induced by a repulsive structure around $11.5$ a.u. makes roaming in acetaldehyde much more probable than in formaldehyde. {\revision Based on the analysis of roaming in formaldehyde we conjecture that the existence of a normally hyperbolic invariant manifold and its stable and unstable manifolds would account for the observed dynamics.}

\section{Conclusion}
\label{sec:concl}
 
 We report a number of roaming trajectories in acetaldehyde with properties that suggest the presence of two distinct roaming pathways for zero total angular momentum.
 For all reported trajectories, CH$_3$ breaks loose from HCO, spends considerable time in a region of nearly constant potential energy and abstracts H from HCO prior to dissociation.
 There are two major features that set these trajectories apart: the maximal distance between CH$_3$ and HCO during roaming, and similarity to in-plane roaming of a highly restricted model of acetaldehyde with 2 degrees of freedom.
 
 As predicted by the restricted model, we find CH$_3$ roaming at approximately  $14.5-22.9$ a.u. from HCO. Additionally we find that phase space structures in the restricted models for acetaldehyde and formaldehyde \cite{mauguiere2015phase} are very similar, suggesting that the roaming mechanism found in formaldehyde is also present in acetaldehyde. This also agrees with the results on mass dependence of roaming in Chesnavich's model \cite{krajnak2018influence}. All of these models consider in-plane roaming.
 
 Unlike predicted by the restricted model, we find CH$_3$ roaming at approximately $9-11.5$ a.u. in the full-dimensional model of acetaldehyde. These trajectories resemble those reported in earlier experimental studies and quasiclassical trajectory calculations. The absence of these trajectories in the restricted in-plane model is in line with the importance of out-of-plane motion for roaming, as emphasised in multiple works \cite{Harding2010roaming, Klippenstein2011roaming, shepler2011roaming}.

 Additionally we find no roaming trajectories that reach a maximal distance between CH$_3$ and HCO between $11.5$ and $14.5$ a.u. While this is not enough to prove that acetaldehyde has two separate roaming mechanisms, the possibility should be certainly considered. Given that formaldehyde was found significantly less likely to dissociate via roaming than acetaldehyde, the presence of an additional roaming mechanism offers a plausible explanation. The evidence presented in this work creates a link between results on roaming in phenomenological models and in realistic molecules, and casts new light on the possibility of in-plane roaming, which is largely considered impossible.

 \section*{Supplementary information}
 Data and videos of roaming trajectories presented in this paper. Each {\itshape hdf5} archive contains an object (numpy.array) called {\itshape traj} that 43 columns of double precision floats: 6 coordinates ($x$, $p_x$, $y$, $p_y$, $z$, $p_z$) corresponding to H, H, H, H, C, C and O respectively and time in the last column. All units are atomic units.

\newpage

\appendix

\section{Restricted acetaldehyde}

\subsection{The Model System}
\label{sec:2dofmodel}

 Derived from a potential for acetaldehyde used in \cite{Han2017} and based on \cite{heazlewood2008roaming} kindly provided by Prof. J. Bowman, Prof. Y.-C. Han and Prof. B. Fu, in this appendix we elaborate on the dynamics of a highly restricted Hamiltonian acetaldehyde system with $2$ degrees of freedom. Our system assumes HCO and CH$_3$ to be rigid bodies with configurations given by the minimum potential energy in the acetaldehyde well. As outcomes of the dynamics of the system we consider the configurations CH$_3$CHO, CH$_3$HCO and the radical products HCO + CH$_3$.
 
 We assume the centre of mass of CH$_3$ to move in a centre of mass frame on the plane defined by HCO. The orientation of CH$_3$ is fixed so that all internal angles as well as angles with respect to the position vector of the centre of mass of CH$_3$ remain constant. In this configuration CH$_3$ is always facing HCO carbon first.

  \begin{figure}
    \centering
    \includegraphics[width=0.4\textwidth]{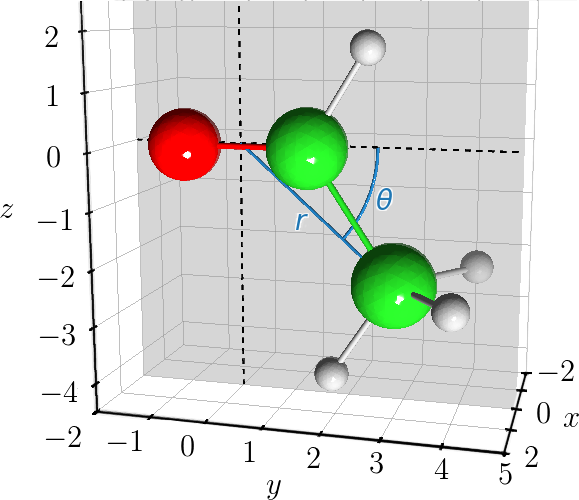}
    \caption{Illustration of polar coordinates describing the position of the centre of mass of CH$_3$ on the HCO plane given by $x=0$.}
    \label{fig:3dconfig}
  \end{figure}
 
 The system is defined by the Hamiltonian {\revision \cite{Ezra2019}}
 \begin{equation}\label{eq:Ham}
  H(r,\theta,p_r, p_\theta) = \frac{p_r^2 }{2 \mu} + \frac{p_\theta^2}{2}\left(\frac{1}{I_{HCO}} + \frac{1}{\mu r^2}\right) + U(r,\theta),
 \end{equation}
 where $r,\theta$ are the polar coordinates of the centre of mass of CH$_3$ on the plane defined by HCO in atomic units and radians respectively (Figure \ref{fig:3dconfig}), $p_r,p_\theta$ are the momenta canonically conjugate to $r, \theta$ respectively, $\mu=18040.89\ m_\text{e}$ is the reduced mass of CH$_3$, $I_{HCO}=81448.95123\ m_\text{e} a_0^2$ is the moment of inertia of HCO and $U$ is the potential energy in hartree derived from the acetaldehyde potential as described above. We investigate the system at a fixed total energy $H=0.01$ a.u. above the dissociation threshold $E_d=-153.440526848047$ a.u.
  
 The Hamiltonian \eqref{eq:Ham} defines the equations of motion
 \begin{equation}\label{eq:Ham_eq}
  \begin{aligned}
  \dot{r} &= \frac{p_r}{\mu},\\
  \dot{p_r} &= \frac{p_\theta^2}{\mu r^3} -\frac{\partial U}{\partial r},\\
  \dot{\theta} &= p_\theta\left(\frac{1}{I_{HCO}} + \frac{1}{\mu r^2}\right),\\
  \dot{p_\theta} &= -\frac{\partial U}{\partial \theta}.
 \end{aligned}
 \end{equation}
 The value of the Hamiltonian is conserved along the solutions and therefore phase space is foliated by energy surfaces of the form $H=E$.
 Total angular momentum is also a conserved quantity, its value in the equations of motion is zero.
 
 \begin{figure}
 \centering
 \includegraphics[width=0.49\textwidth]{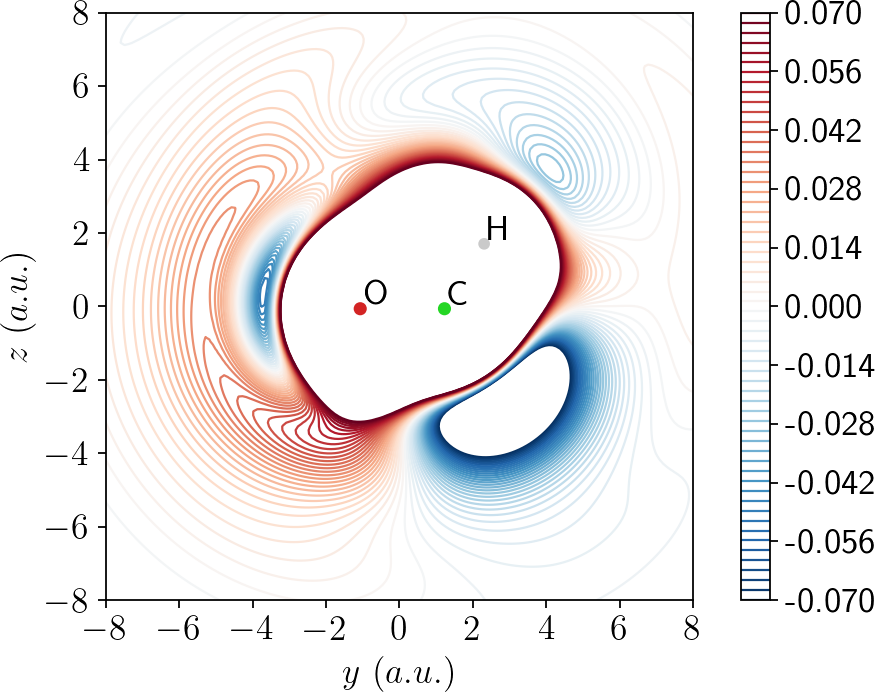}
 \caption{Contour plot of the potential energy surface showing potential wells (blue) and areas of high potential (red).}
 \label{fig:pot}
 \end{figure}
  
\begin{figure}
 \centering
 \includegraphics[width=0.51\textwidth]{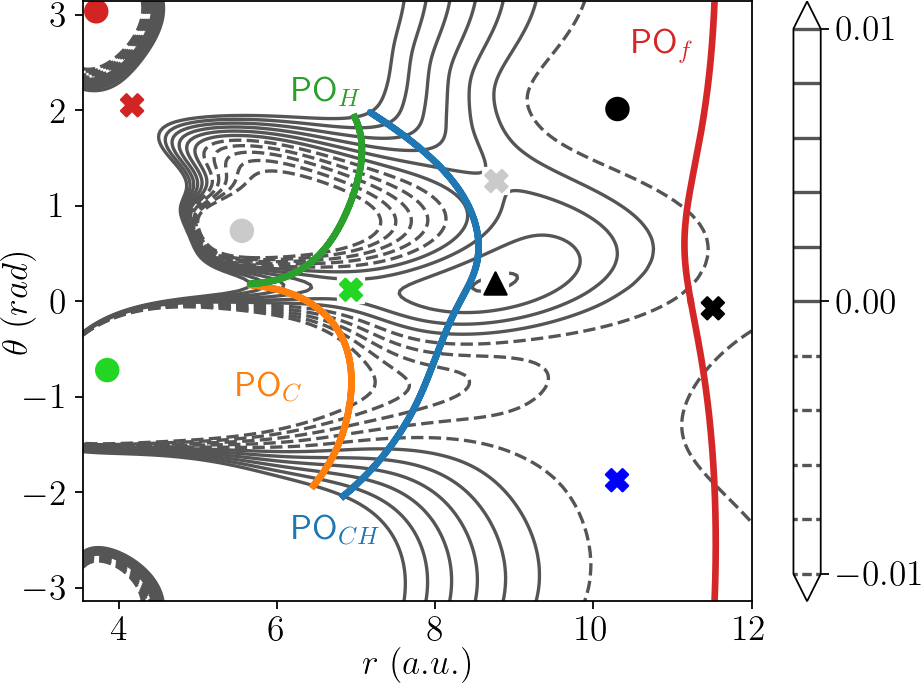}
 \caption{Configuration space projections of relevant periodic orbits PO$_C$, PO$_H$, PO$_{CH}$, PO$_f$, with the outer periodic orbits PO$_O$ ($r\approx22.9$) left out for the sake of readability. Critical points of the potential are shown and listed in Tab. \ref{tab:critpoints}. Local minima are marked with a circle, index-1 saddles with a cross and an index-2 saddle with a triangle.}
 \label{fig:orbits}
 \end{figure}

 A contour plot of $U$ in cartesian coordinates is shown in Figure \ref{fig:pot}. Critical points, their coordinates and energies (with respect to $E_d$) are listed in Tab. \ref{tab:critpoints} and displayed in Figure \ref{fig:orbits}.

  \begin{table}
    \begin{center}
    \begin{tabular}{c|c|c|c|c}
    Energy (hartree) & $r$ (a.u.) & $\theta$ (radians) & Significance & Label in Figure \ref{fig:orbits}\\
    \hline
    $-0.16116$ & $3.85430$ & $-0.72053$ & C-well (acetaldehyde) & {\color{C} $\medblackcircle$}\\
    $-0.02910$ & $5.55601$ & $\Ph0.73850$ & H-well & {\color{H} $\medblackcircle$} \\
    $\Ph0.00300$ & $6.93530$  & $\Ph0.12084$ & saddle between C-well and H-well & {\color{C} \ding{54}}\\
    $\Ph0.00188$ & $8.77285$ & $\Ph1.26013$ & saddle between H-well and flat region&  {\color{H} \ding{54}}\\
    $\Ph0.00616$ & $8.75635$ & $\Ph0.19413$ & index-2 saddle & $\medblacktriangleup$ \\
    
    \hline
    $-0.00381$ & $10.30339$ &  $\Ph2.01440$ & well in the flat region& $\medblackcircle$ \\
    $-0.00246$ & $10.29796$ & $-1.87346$ & saddle in the flat region&  {\color{blue} \ding{54}}\\
    $-0.00233$ & $11.50929$ & $-0.07327$ & saddle in the flat region& \ding{54}\\     
    $-0.05092$ & $3.71536$ & $\Ph3.03771$ & O-well & {\color{O} $\medblackcircle$} \\
    $\Ph0.01752$ & $\phantom{1}4.16820$ & $\Ph2.05342$ & saddle between H-well and O-well& {\color{O} \ding{54}} \\
    \end{tabular}
    \end{center}
    \caption{\label{tab:critpoints} Critical points of the potential $U(r, \theta)$. Positions are shown in Figure \ref{fig:orbits}}
  \end{table}
   
 There are three potential wells near HCO, the deepest of which is the acetaldehyde well near C (further {\revision referred} to as the C-well). The potential well via which acetaldehyde dissociates into molecular products near H (H-well), also relevant to roaming. These wells are separated by an index-1 saddle and an index-2 saddle. Both wells are connected to the flat region via index-1 saddles with energies below the dissociation threshold. A third potential well near O is inaccessible at $H=0.01$ a.u.  
  
\begin{figure}
 \centering
 \includegraphics[width=\textwidth]{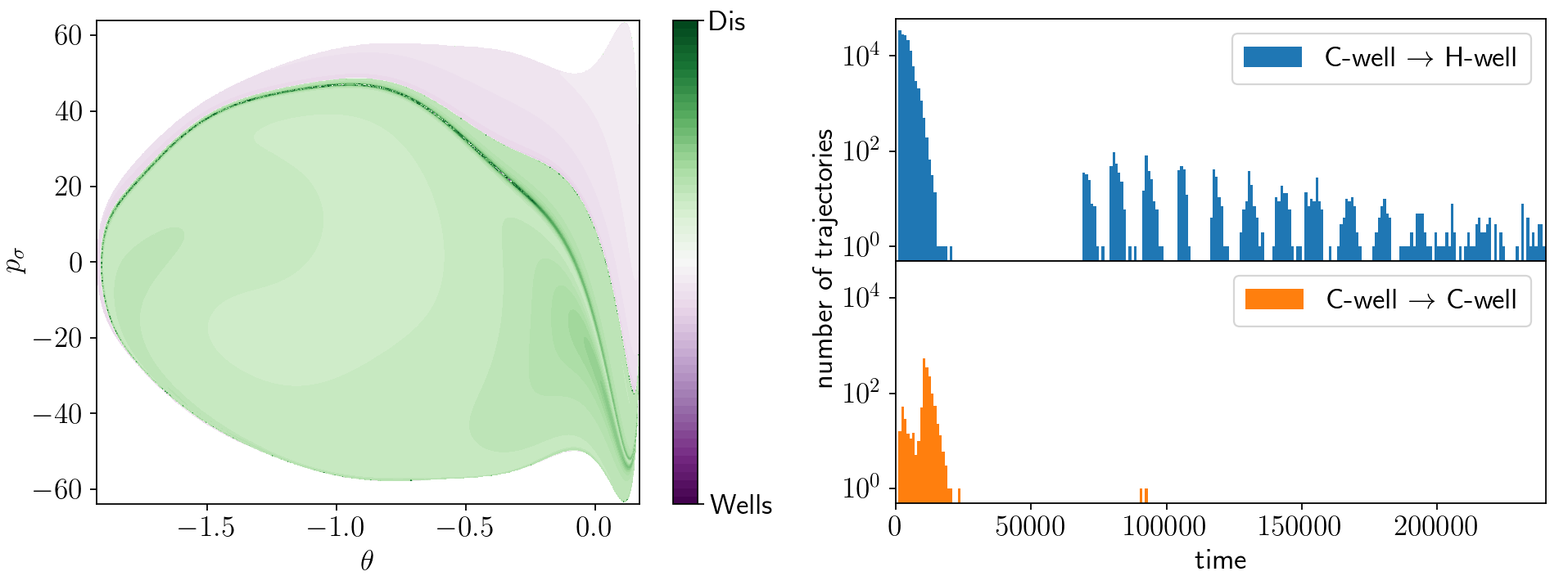}
 \caption{Escape/capture times of trajectories leaving the C-well initialised uniformly on the outward hemisphere of DS$_C$ for $H=0.01$. Left: Initial conditions that lead to radical dissociation are shown in shades of green, trajectories that enter either of the wells in purple. Dark shades indicate long escape/capture time, light shades short. $p_\sigma$ on the vertical axis is the momentum canonically conjugate to $\theta$ on DS$_C$. Right: Histograms of capture times of the corresponding trajectories terminated in the H-well and the C-well respectively. The vertical axis is in log scale.}
 \label{fig:escape_C}
 \end{figure}

 Trajectory simulations suggest, and a phase space analysis confirms, the presence of two different mechanisms via which the restricted system dissociates into molecular products - one passing near the potential saddle point and one that initially leads CH$_3$ into the flat region. Figure \ref{fig:escape_C} shows the escape/capture times of $767266$ trajectories leaving the C-well initialised uniformly on the appropriate dividing surface defined in the Dividing surfaces section. Of those $627506$ dissociated, $1522$ entered the C-well, $138238$ entered the H-well - but only $1221$ took longer than $50000$ time units. We explain how stable and unstable invariant manifolds of unstable periodic orbits convey trajectories between the two wells in phase space. 
 
\subsection{Relevant periodic orbits}
\label{sec:identify}

 We identified several vibrational and rotational periodic orbits that govern the dynamics of the restricted acetaldehyde, the configuration space projection of four of which are shown in Figure \ref{fig:orbits}. Note that rotational periodic orbits always come in pairs with identical configuration space projection and opposite orientation - one rotating clockwise and one counter-clockwise. Most of these orbits and their dynamical significance are similar to those found in other roaming systems such as formaldehyde \cite{mauguiere2015phase}, Chesnavich's CH$_4^+$ model \cite{mauguiere2014multiple,krajnak2018phase} or the double Morse system \cite{montoya2020morse}. 
 These are
 \begin{itemize}
  \item a pair of nearly circular ($r\approx22.9$) unstable outer periodic orbits which marks the point of radical dissociation,
  \item an unstable orbit delimiting the C-well,
  \item an unstable orbit delimiting the H-well,
  \item a pair of periodic orbits in the flat region.
 \end{itemize}
 
 We refer to these orbits as the outer orbits (PO$_O$), PO$_C$, PO$_H$ and PO$_f$ respectively. None of these orbits is related to a saddle point of the potential in an obvious way. The PO$_O$ pair is due to a centrifugal barrier \cite{krajnak2018phase}.
 
 The dividing surfaces DS$_C$, DS$_H$ and DS$_f$ are defined as the surfaces on the $H=0.01$ energy surface that have the same configuration space projections as PO$_C$, PO$_H$ and PO$_f$ respectively. This construction is due to \cite{PollakPechukas78,pechukas1979} and details are provided in the dividing surfaces section. DS$_C$ and DS$_H$ delimit the analogue of the respective potential wells in phase space. We refrain from using the term 'transition state' due to its ambiguous use in the literature for a potential saddle point, the associated periodic orbit (which depending on energy may or may not be unstable) and the associated dividing surface.
 
 Unlike formaldehyde, Chesnavich's model and the double Morse, a crucial role in restricted acetaldehyde is played by another orbit, referred to as PO$_{CH}$:
  \begin{itemize}
  \item an unstable periodic orbit that lies between the C- and H-wells and the flat region.
 \end{itemize}
 Unstable periodic orbits are located in phase space bottlenecks, thanks to which the corresponding DS has locally minimal flux \cite{pechukas1979,waalkens2004}. This bottleneck allows certain trajectories to pass, while preventing others from doing so. The selection is due to the stable and unstable invariant manifolds \cite{Rom-Kedar90,Uzer02,Waalkens04} of the unstable periodic orbit. This leads to two distinct mechanism of transport between the C-well and the H-well: a direct mechanism that takes trajectories over a potential barrier near the potential saddle point from one well almost immediately to the other and an indirect mechanism that leads trajectories via the flat region corresponding to roaming.
 
 Stable and unstable invariant manifolds of unstable periodic orbits are calculated using a monodromy matrix as detailed in \cite{krajnak2020hydrogen}.
 
\subsection{Roaming in Restricted Acetaldehyde}
\label{sec:rra}

 Each DS can be crossed by a trajectory in two directions, roughly speaking, in the direction of increasing or decreasing radius. These two halves are separated by the respective periodic orbits and will be referred to as inward and outward respectively. More details are in the Dividing surfaces section.
 
 Roaming in the restricted acetaldehyde molecule is the transport of points from the outward half of DS$_C$ to the inward half of DS$_H$ via both halves of DS$_f$. This definition is consistent with the definition of roaming in formaldehyde \cite{mauguiere2015phase}.
 
 Unstable periodic orbits have cylindrical stable and unstable invariant manifolds, that act as barriers to dynamics and mediate passage through the corresponding phase space bottleneck. The relevance of unstable periodic orbits was recognised by Pechukas and Pollak \cite{PollakPechukas78,pechukas1979} in the context of transition state theory. Dynamics in higher dimensional systems is governed by the manifolds of normally hyperbolic invariant manifolds \cite{Fenichel1971,Wiggins94}, of which an unstable periodic orbit is a $1$-dimensional example.
 
 A stable invariant manifold $\mathcal{S}^\gamma$ of an unstable periodic orbit $\gamma$ consist of trajectories that converge to $\gamma$ as $t\rightarrow\infty$. Due to its cylindrical structure, $\mathcal{S}^\gamma$ encloses a volume of trajectories that is conveyed to DS$_\gamma$ at the narrowest part of the bottleneck. Trajectories not enclosed by $\mathcal{S}^\gamma$ do not reach DS$_\gamma$. An unstable invariant manifold $\mathcal{U}^\gamma$ consist of trajectories that converge to $\gamma$ as $t\rightarrow-\infty$, but diverge in forward time. $\mathcal{U}^\gamma$ encloses all trajectories that crossed DS$_\gamma$ and guide them away. The roles of $\mathcal{S}^\gamma$ and $\mathcal{U}^\gamma$ reverse in backward time.

 Each invariant manifold consists of two branches: one conveying trajectories to/from the inward half of DS$_\gamma$ and one to/from the outward half.
 By definition a stable invariant manifold cannot intersect itself or another stable invariant manifold and the analogue is true for unstable invariant manifolds.
 
 On the other hand, the intersection of $\mathcal{U}^{C}$ and $\mathcal{S}^H$ will guide trajectories from the outward half of DS$_C$ to the inward half of DS$_H$, which is crucial for roaming.
 Radical dissociation is governed by the intersection of $\mathcal{U}^C$ and $\mathcal{S}^O$.
   
 \begin{figure}
 \centering
  \includegraphics[width=\textwidth]{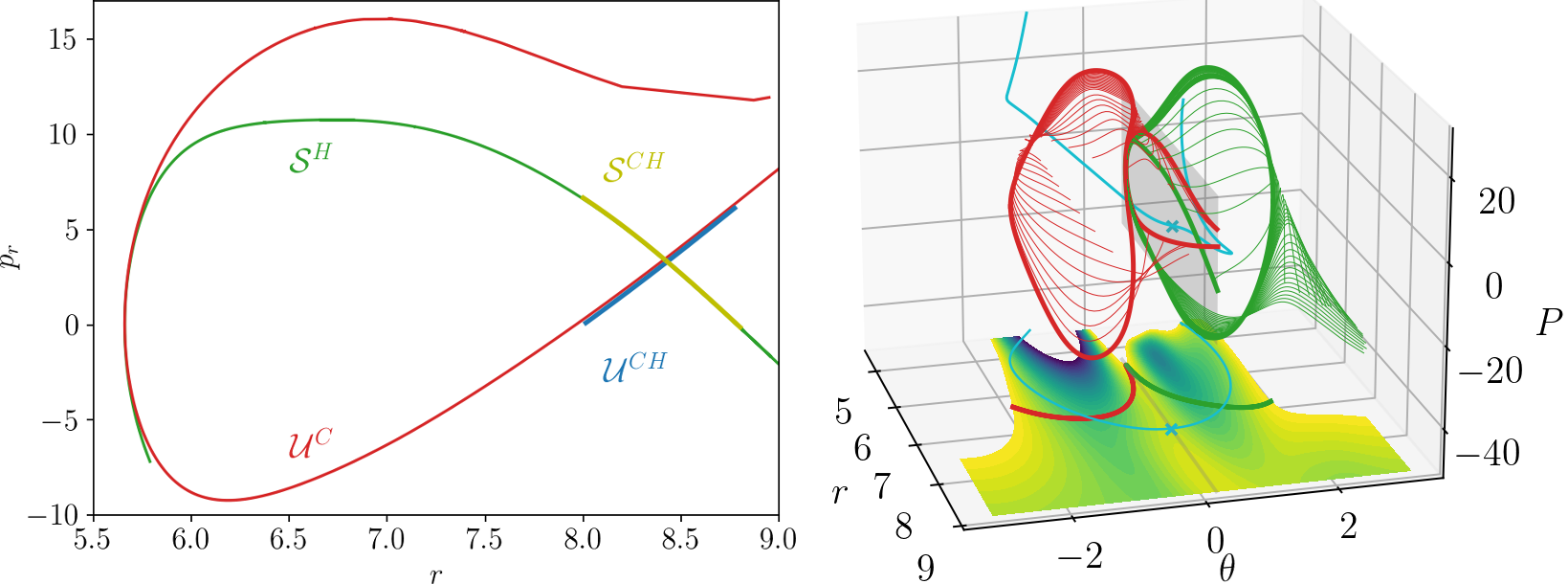}
 \caption{The mechanism conveying trajectories from the C-well to the H-Well over the potential barrier. Left: Intersection of $\mathcal{U}^{C}$ and $\mathcal{S}^H$ on $\theta=0.177$, $p_\theta>0$. Right: A $3$-dimensional representation of the intersection of $\mathcal{U}^{C}$ and $\mathcal{S}^H$. For visualisation purposes, values on the vertical axis correspond to a $(r, \theta)$-dependent linear combination $P$ of $p_r$ and $p_\theta$, such that $P$=$p_r$ on $\theta=0.177$. The section the invariant manifold cylinders on $\theta=0.177$, $p_\theta>0$ is highlighted. A representative trajectory is shown in cyan and the intersection with $\theta=0.177$, $p_\theta>0$ is marked with a cross. Configuration space projections of the periodic orbits and trajectory on a contour plot of the potential energy surface are added for illustration.}
 \label{fig:manifs_th0177}
 \end{figure}
 
 \begin{figure}
 \centering
 \includegraphics[width=\textwidth]{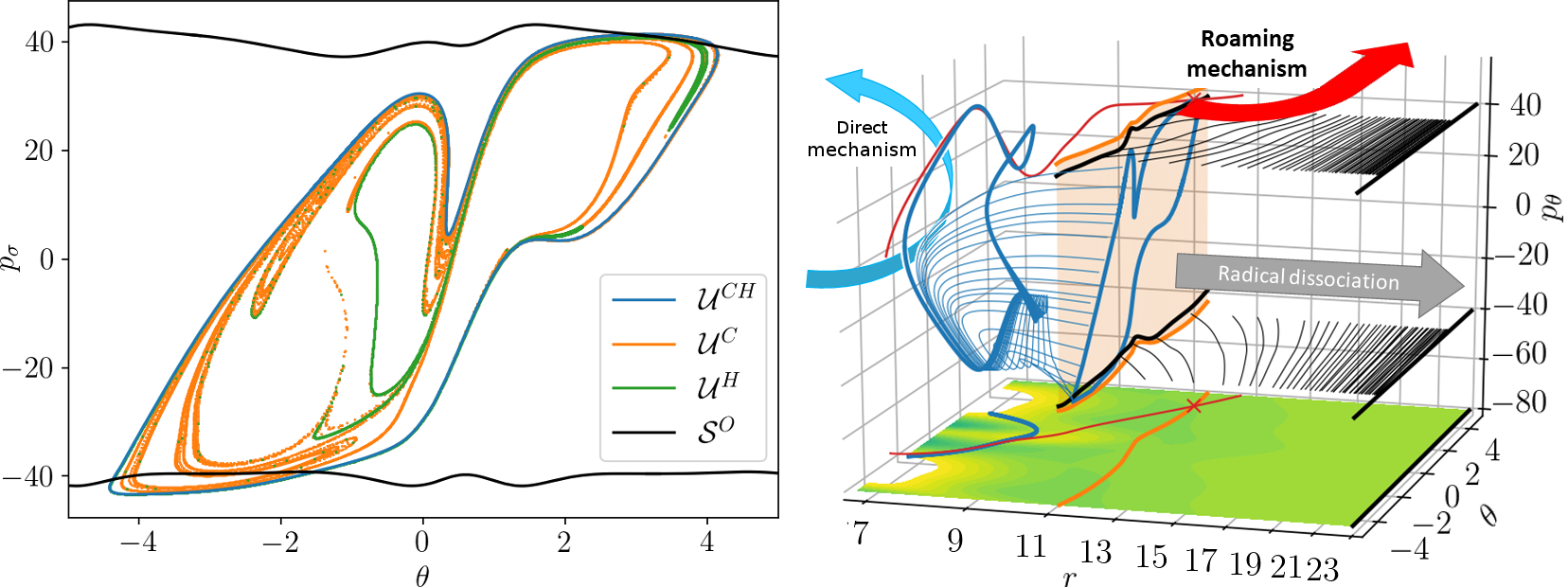}
 \caption{The roaming mechanism conveying trajectories from the C-well to the H-Well via the flat region. Left: Intersection of $\mathcal{U}^{CH}$ and $\mathcal{U}^{C}$ with $\mathcal{S}^O$ on the outward annulus of DS$_f$. Right: A $3$-dimensional projection of the intersection of $\mathcal{U}^{CH}$ and $\mathcal{S}^O$. The section the invariant manifold cylinders on the outward annulus of DS$_f$ is highlighted. A representative trajectory is shown in red and the intersection with the outward annulus of DS$_f$ is marked with a cross. Configuration space projections of the periodic orbits and trajectory on a contour plot of the potential energy surface are added for illustration. The $r$-axis is in log scale for better readability.}
 \label{fig:manifs_A}
 \end{figure}

 The PO$_{CH}$ sets restricted acetaldehyde apart from the likes of formaldehyde, Chesnavich's model or the double Morse. It allows $\mathcal{U}^C$ and $\mathcal{S}^H$ to intersect in two distinct regions:
 \begin{itemize}
  \item the region close to the wells near DS$_C$ and DS$_H$,
  \item the flat region beyond PO$_{CH}$,
 \end{itemize}
 separated by a phase space bottleneck associated with PO$_{CH}$. These two separate intersections form two different mechanisms for transport between the outward half of DS$_C$ and the inward half of DS$_H$:
 \begin{itemize}
  \item a direct mechanism that takes trajectories from DS$_C$ almost immediately to DS$_H$,
  \item an indirect mechanism that leads trajectories first to DS$_f$ in the flat region - roaming.
 \end{itemize}
 
 Figure \ref{fig:manifs_th0177} shows the intersection of $\mathcal{U}^C$ and $\mathcal{S}^H$ on the surface of section $\theta=0.177$, $p_\theta>0$. These manifold intersections should be interpreted as reactive islands that do not form closed circles, because some trajectories on the manifolds are tangent to $\theta=0.177$, $p_\theta>0$. $\mathcal{U}^{CH}$ and $\mathcal{S}^{CH}$ enclose the intersection of $\mathcal{U}^C$ and $\mathcal{S}^H$, highlighting the fact that it does not pass through the phase space bottleneck associated with PO$_{CH}$. Only trajectories above $\mathcal{S}^{CH}$ and below $\mathcal{U}^{CH}$ pass through the bottleneck. A 3-dimensional representation of the mechanism is shown in Figure \ref{fig:manifs_th0177}.
 
 Evidence for the indirect, roaming mechanism, can be found in the flat region. As shown in the analysis of Chesnavich's CH$_4^+$ system \cite{krajnak2018phase, krajnak2018influence}, while the existence of a centrifugal barrier is necessary for roaming, the mechanism for roaming relies on the intersection of $\mathcal{U}^C$ and $\mathcal{S}^O$. Figure \ref{fig:manifs_A} shows the intersection of $\mathcal{U}^{C}$ and $\mathcal{S}^O$ on the outward annulus of DS$_f$ and a $3$-dimensional projection putting the intersection into wider mechanistic context.
 
 The mechanism works in a simple way - all trajectories contained in the interior of both $\mathcal{U}^C$ and $\mathcal{S}^O$ dissociate immediately. Due to their intersection, there is a part of the interior of $\mathcal{U}^C$ that does not dissociate, but stays in the flat region. It exhibits chaotic behaviour due to the presence of pieces of stable and unstable invariant manifolds of all unstable periodic orbits that cause severe distortion. Trajectories contained in this part of $\mathcal{U}^C$ may roam in the traditional sense or exhibit transient roaming \cite{Endo2020}.
 
 The passage from the flat region to DS$_{H}$ and the H-well is mediated by $\mathcal{S}^H$. Due to time reversal symmetry, $\mathcal{U}^H$ and $\mathcal{S}^O$ intersect if and only if $\mathcal{S}^H$ and $\mathcal{U}^O$ intersect. The intersection of $\mathcal{U}^H$ and $\mathcal{S}^O$ is shown in Figure \ref{fig:manifs_A}.
 
 We remark that a similar argument can be used to deduce the mechanism responsible for the passage from the flat region to DS$_{C}$ and the C-well, that explains the two orange outliers in the histogram in Figure \ref{fig:escape_C}.
 
 Animations of representative trajectories for the direct mechanism, the roaming mechanism and radical dissociation are provided in Supplementary Information. 
 
 The histogram in Figure \ref{fig:escape_C} further shows that roaming trajectories enter the H-well in `packets' rather than being uniformly randomised. This agrees with a lobe dynamics \cite{Rom-Kedar90} point of view, where trajectories contained in a lobe (`packet') share dynamical properties, such as in this case the number of rotations of CH$_3$ around HCO. The peaks of the roaming time distribution are approximately $12616$ time units apart, which corresponds to the period of PO$_O$. This also agrees with the existing classification of dynamical behaviour at constant kinetic energy \cite{krajnak2019isokinetic}. We remark that the maximal time in the histogram is an artificial bound introduced for numerical integration. Lobe dynamics teaches us that no maximum exists due to the fractal structure of invariant manifold intersections.
  
 \subsection{Finding periodic orbits for Restricted Acetaldehyde}
 \label{sec:finding}
 
 We find vibrational and rotational periodic orbits in different ways. Vibrational periodic orbits are characterised by two points at which kinetic energy vanishes and the all energy is in the potential component. These orbits can be found via a bisection procedure along the equipotential line(s) $U(r,\theta)=E$ initiated with $p_r=p_\theta=0$.
 
 Rotational periodic orbits require a different approach. Due to the absence of symmetries of this system, we used Lagrangian descriptors \cite{madrid2009ld} in the sense as \cite{krajnak2020manifld} to identify locations of nearly $\theta$-independent dynamics and subsequently employed the binary contraction method \cite{bardakcioglu2018binary} to locate PO$_O$ and PO$_f$.
 
\begin{figure}
 \centering
 \includegraphics[width=0.49\textwidth]{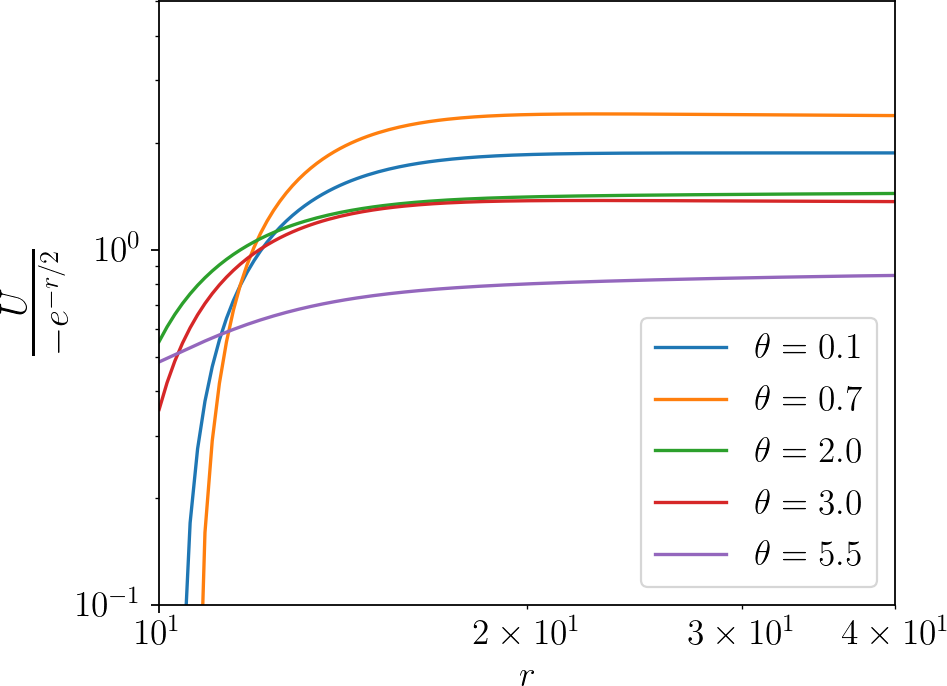}
 \caption{Sections of $\frac{U}{-e^{-r/2}}$ for constant values of the angle $\theta$ showing that $U\in O(e^{-r/2})$ as $r\rightarrow\infty$. Both axes are in logarithmic scale. A rate of convergence of the potential to the dissociation energy faster than $r^{-2}$ assures the existence of a centrifugal barrier.}
 \label{fig:pot_sec}
 \end{figure}

 We remark that the existence of PO$_O$ is due to a centrifugal barrier that exists in any system for sufficiently low energies, provided the potential $V$ approaches the dissociation energy at least $V\in o(r^{-2})$ as $r\rightarrow\infty$  \cite{krajnak2018phase}. In fact $U$ converges faster, $U\in O(e^{-r/2})$ as shown in Figure \ref{fig:pot_sec}.
  
 \subsection{Construction Dividing Surfaces for Restricted Acetaldehyde}
\label{sec:DS}
 
 Dividing surfaces (DS) are constructed using periodic orbits as follows. We remark that while the instability of a periodic orbit assures locally minimal flux across the DS, the construction can be done using any periodic orbit.
 
 Let $\Gamma$ be a periodic orbit. For any point $(r_\Gamma, \theta_\Gamma)$ on the configuration space projection of $\Gamma$, consider all values of $p_r$ and $p_\theta$, such that the point $(r_\Gamma, p_r, \theta_\Gamma, p_\theta)$ lies on the energy surface. All such points project onto $(r_\Gamma, \theta_\Gamma)$ and satisfy
 \begin{equation*}
  \frac{p_r^2 }{2 \mu} + \frac{p_\theta^2}{2}\left(\frac{1}{I_{HCO}} + \frac{1}{\mu r_\Gamma^2}\right) = H - U(r_\Gamma,\theta_\Gamma).
 \end{equation*}
 On the energy surface this set is a circle and collapses to a point if $H = U(r_\Gamma,\theta_\Gamma)$. We define the DS as the union of these sets over all points $(r_\Gamma, \theta_\Gamma)$.
  
 As $H = U(r,\theta)$ at the turning points of the corresponding orbits, both DS$_C$ and DS$_H$ are spheres consisting of a continuum of circles and two points at the poles. The construction above is known to provide a DS that is transversal to the Hamiltonian vector field, in our case \eqref{eq:Ham_eq}, except at the periodic orbit \cite{pechukas1979,waalkens2004}.
  
 The construction of a DS is identical for a pair of rotational periodic orbits, such as PO$_f$ or PO$_O$, since the orbits share a configuration space projection and only differ in the sign of their momenta. Since  $H > U(r,\theta)$ along rotational periodic orbits, the DS is known to be a torus and transversal to the Hamiltonian vector field  \cite{mauguiere2016ozone}.
 
 Note that DSs are not unique, as any deformation of the above-constructed DS that does not violate transversality to the vector field is also a DS.
 
 Each of the above-defined DSs can be split in two halves: spheres are split into two hemispheres bounded by their respective periodic orbit, tori into two annuli bounded by their two orbits. Due to transversality, each half must be crossed by trajectories in the same direction, for example in the case of DS$_C$ escape from the bound states of restricted acetaldehyde. Due to conservation of volume, the flux across the two halves must be equal in magnitude and opposite in direction.
 
 In the literature, the hemisphere of DS$_C$ through which trajectories escape from the bound states of restricted acetaldehyde would be referred to as the outward hemisphere. Similarly the inward hemisphere is the one through which trajectories enter the bound states of restricted acetaldehyde.

\subsection{Conclusion for the restricted system}
\label{sec:concl_res}

 We found two mechanisms consisting of phase space structures that lead a restricted acetaldehyde molecule to molecular products - passage over a potential barrier and roaming. 
 The $2$ degree-of-freedom system considered in this appendix does not capture roaming in its entirety, that is with out-of-plane motion of the CH$_3$ radical or the process of hydrogen abstraction. Describing these processes requires a higher dimensional system and correction terms\cite{Harding2010roaming}. A complete understanding of the mechanisms governing the dynamics of restricted subsystems is crucial for future analyses of the dynamics of higher dimensional systems.
 
 This work considers the dynamics at a single total energy, $0.01$ a.u. above the dissociation threshold. At this energy we find passage over the potential barrier to dominate the roaming mechanism, both eclipsed by radical dissociation, similar to the proportions reported for formaldehyde \cite{mauguiere2015phase}. 
 Qualitatively we find that the roaming mechanisms of restricted formaldehyde and restricted acetaldehyde bear remarkable similarities, the only differences induced by PO$_{CH}$ and associated bottleneck. The PO$_{CH}$ in restricted acetaldehyde, which is absent in formaldehyde, allows for a direct transport passage between the wells.
 
 Results on persistence of normally hyperbolic invariant manifolds under perturbations \cite{Fenichel1971,Wiggins94} assure that all of the unstable periodic orbits involved in roaming persist for a range of energies around $0.01$ a.u. Similarly our results remain valid under small variations of the underlying potential energy surface, small changes to the geometry of the rigid bodies or when the system is weakly coupled to harmonic bath modes \cite{Naik2020systembath}.
 
 As long as the periodic orbits creating the template for reaction mechanisms do not change, the reaction mechanisms themselves do not change qualitatively. For larger perturbations bifurcations of periodic orbits are possible with resulting changes to reaction mechanisms. The connection between bifurcations and changes in reaction mechanisms is the topic of ongoing research.
  
 The existence of a bottleneck separating the wells from the flat region in restricted acetaldehyde gives rise to a new reaction mechanism transporting trajectories between the wells. In the case of more complicated potential energy surfaces, every additional bottleneck can be expected to give rise to at least one additional reaction mechanism. The necessary passage through multiple bottlenecks may decrease the proportion of roaming among competing mechanisms for molecular dissociation or even give rise to multiple roaming mechanisms, but as long as the equivalents of $\mathcal{U}^{C}$ and $\mathcal{S}^O$ intersect, roaming will be present.

\section*{Data availability}
  The data that support the finding of this study are available in the Supplementary material.

\section*{Acknowledgments}
 
  We thank Prof. J. Bowman, Prof. Y.-C. Han and Prof. B. Fu for providing us with the potential energy surface of acetaldehyde. We are grateful to Prof. Gregory S. Ezra for his comments and insights on roaming, and chemical dynamics in general.
  We acknowledge the support of  EPSRC Grant no. EP/P021123/1 and Office of Naval Research (Grant No.~N000141712220).
  We also acknowledge the high performance computing cluster Cream at the School of Mathematics and BlueCrystal at the Advanced Computing Research Centre of University of Bristol.

\section*{Author contributions statement}
V.K. performed the calculations and analysis, V.K. and S.W. discussed the results and reviewed the manuscript.

\section*{Additional information}
\textbf{Competing interests} The authors declare no competing interests.

\bibliography{acetaldehyde}
\end{document}